\documentclass[12pt,draft]{article}
\usepackage[cp1250]{inputenc}
\usepackage{amssymb}
\usepackage{amsmath}
\usepackage{fancyhdr}
\usepackage{amsthm}

\font\teneufm=eufm10 \font\seveneufm=eufm7 \font\fiveeufm=eufm5
\newfam\eufmfam
\textfont\eufmfam=\teneufm \scriptfont\eufmfam=\seveneufm
\scriptscriptfont\eufmfam=\fiveeufm

\def\bee{\begin{eqnarray}}
\def\bes{\begin{eqnarray*}}
\def\eee{\end{eqnarray}}
\def\ees{\end{eqnarray*}}

\begin{document}

\def\O{{\hbox{\boldmath $O$}}}
\def\R{{\hbox{\boldmath $R$}}}
\def\C{{\hbox{\boldmath $C$}}}

\def\le{{\langle}}
\def\re{{\rangle}}

\title{Simple Moufang loops and alternative algebras}

\author{Sandu N. I.}
\date{}

\maketitle
\begin{abstract}
Let a Moufang loop $Q$ contain a non-unitary subloop, which is a
simple loop. Then $Q$  is not embedded into a loop of invertible
elements of any alternative algebra.
\end{abstract}

\noindent {\bf Mathematics Subject Classification (1991):}  17D05,
20N05.

\noindent {\bf Key words:} simple Moufang loop,  alternative
algebra.
\bigskip

The alternative algebras satisfy the Moufang identities that
define the Moufang loops and the set $U(A)$ of all invertible
elements of any alternative algebra $A$ forms a Moufang loop with
respect to multiplication \cite{M}. Therefore,  the question
raised in \cite{G} is natural.
\smallskip

{\bf Question.\ }{\em Is it true that any Moufang loop can be
imbedded into a loop of type $U(A)$ for a suitable unital
alternative algebra $A$?}
\smallskip

In general, the answer to this question is negative. In \cite{Sh}
the Moufang loop $U(\O)/\R^{\star}$ for the algebra $\O$ of
classical Cayley algebra (the Cayley-Dikson algebra $\C(-1, -1,
-1)$ in the notation of \cite{ZSSS}) over the real field $\R$ and
the similar loop for the Cayley-Dickson algebra  over the finite
field $GF(p^2),\ p>2,$ are not imbeddable into the loops of type
$U(A)$.
\smallskip

The field $CF(p^2)$ is split. Any quadratic equation with
coefficients from $CF(p)$ is solvable in $CF(p^2)$. Then the field
$CF(p^2)$ is closed under root operation and by [5, pag. 475,
Theorem] the corresponding Moufang loop for $CF(p^2)$ is simple.
\smallskip

Now we prove that the Moufang loop $U(O)/\R^{\star}$ also is
simple for classical Cayley algebra $\O$. Really, suppose
contrary, that the loop $U(\O)/\R^{\star}$ is non-simple. Let
$\overline{H}$ be a proper normal subloop of loop
$U(O)/\R^{\star}$ and let $H$ be the inverse image of
$\overline{H}$ under homomorphism loops $U(\O) \rightarrow
U(\O)/\overline{H}$. Then $H$ will be a proper normal subloop of
$U(\O)$.

The Cayley algebra $\O$ is a $8$-dimensional corp \cite{K}. Then
$U(\O) = \O \backslash \{0\}$. Let $1 = e_0, e_1, \ldots, e_7$ be
the canonical basis of $\O$. The normal subloop $H$ of $U(\O)$ is
proper. Then $e_j \notin H$ for some $e_j \in \{e_0, \ldots,
e_7\}$.

Any element in $\O$ has a form $\sum_{i=0}^7\alpha_ie_i$, where
$\alpha_i \in \R$. As $U(\O) = \O \backslash \{0\}$ then the
mapping $\varphi: \sum_{i=0}^7\alpha_ie_i \rightarrow
\sum_{i=0}^7\alpha_ie_iH$ is a homomorphism of algebra $\O$.
Obviously, $H = U(\O) \cap (1 + \text{ker}$ $\varphi)$. Then from
$e_j \notin H$ it follows that $e_j - 1 \notin \text{ker}$
$\varphi$. Hence $\text{ker}$ $\varphi$ is a proper ideal of $\O$,
i.e. $\O$ is a non-simple alternative algebra. But by Kleinfeld
Theorem the Cayley-Dickson algebras  and only they can be simple
alternative algebras \cite{ZSSS}. We get a contradiction.
Consequently, the Moufang loop $U(O)/R^{\star}$ is simple.
\smallskip

The following Theorem together with those proved above generalize
the main  result from \cite{Sh}.  Let us prove first the
following:
\smallskip

\textbf{Lemma.} \textit{Let $I$ be an ideal of alternative algebra
$A$ with unit $1$ and let $H$ be a subloop of Moufang loop $U(A)$.
Then $K = H \cap (1 + I)$ is a normal subloop of loop $H$.}
\smallskip

\textbf{Proof.} The homomorphism of algebras $A \rightarrow A/I$
induce a homomorphism $\varphi$ of loop $H$. Denote $\varphi(H,
\cdot) = (\overline H,\star)$. Any Moufang loop is a $IP$-loop,
i.e. satisfies the identities $x^{-1}\cdot xy = y$, $yx\cdot
x^{-1} = y$, where $x^{-1}x = xx^{-1} = 1$. From identity
$x^{-1}\cdot xy = y$ it follows that $\varphi(x^{-1}) = (\varphi
x)^{-1}$ and $(\varphi x^{-1})\star (\varphi x\star\varphi y) =
\varphi y, (\varphi x)^{-1}\star(\varphi x\star\varphi y) =
\varphi y, \overline x^{-1}\star(\overline x\star \overline y) =
\overline y$.

Let $\overline a, \overline b \in \overline H$. It is obvious that
the equation $\overline a\star x = \overline b$ is always solvable
and as $\overline a^{-1}\star(\overline a\star x) = \overline
a^{-1}\star \overline b, x = \overline a^{-1}\star \overline b$,
then it is uniquely solvable. It can be shown by analogy that the
equation $y\star \overline a = \overline b$ is also uniquely
solvable. Hence $(\overline H,\star)$  is a loop. Then
$\text{ker}$ $\varphi$ $= H \cap (1 + I)$ is a normal subloop of
loop $H$, as required.
\smallskip

\textbf{Theorem}. \textit{Let a Moufang loop $Q$ contain a
non-unitary subloop, which is a simple loop.  Then the loop $Q$ is
not imbedded into the loop of type $\mathcal{U}(A)$ for a suitable
unital alternative $F$-algebra $A$, where $F$ is an associative
commutative ring with unit.}
\smallskip

\textbf{Proof.} Obviously, it is sufficient to consider that $Q$
is a simple Moufang loop. In process of proof will be used without
reference some definitions and results from theory of alternative
algebras from \cite{ZSSS}.
\smallskip

Assume that $Q$ is imbedded into a loop $\mathcal{U}(A)$ for a
certain alternative algebra $A$. We will identify the elements
from $Q$ with their images in $Q$.

Let $F\{Q\}$ be the submodule of $F$-module $A$ generated by set
$\{g \vert g \in Q\}$. Finite sums $\sum_{g \in Q}\alpha_gg$,
where $\alpha_g \in F$, and only they are elements in $F\{Q\}$.
Obviously, $F\{Q\}$ is a subalgebra of algebra $A$, $Q$ is a
subloop of loop $U(F\{Q\})$. If $I$ is a proper ideal of $F\{Q\}$
then $g \notin I$ for some $g \in Q$. In this case by Lemma the
ideal $I$ of $F\{Q\}$ induces the proper normal subloop $K = Q
\cap (1 + I)$ of loop $Q$.

Let $I_1, I_2$ be a proper ideals of algebra $F\{Q\}$ and let
$K_1, K_2$  be the proper normal subloops of loop $Q$
corresponding to ideals $I_1, I_2$ by Lemma. The sum $I_1 + I_2$
is the minimal ideal of $F\{Q\}$ containing the ideals $I_1, I_2$,
the product $K_1K_2$ is the minimal normal subloop of $Q$
containing the normal subloops $K_1, K_2$. Hence the ideal $I_1 +
I_2$ induces the normal subloop $K_1K_2$ of loop $Q$.

The loop $Q$ is simple. Then $K_1 = K_2 = 1$. Hence the sum $I_1 +
I_2$ of any proper ideals $I_1, I_2$ of $F\{Q\}$ induces in loop
$Q$ the unitary subloop $1$, i.e. induces the identical mapping in
$Q$.

It is easy to see that the sum $S$ of all proper ideals of
$F\{Q\}$ induces in loop $Q$ the identical mapping in $Q$. Hence
the algebra $F(Q) = F\{Q\}/S$ is non-trivial, simple and $Q$ is a
subloop of loop $U(F(Q))$. The elements $g \in Q$ are invertible,
then the elements $1 - g$ are quasiregular in algebra $F(Q)$.
Hence the Smiley radical $\mathcal{S}(F(Q))$ of algebra $F(Q)$ is
non-trivial, $\mathcal{S}(F(Q)) \neq (0)$. Remind that the Smiley
radical $\mathcal{S}(A)$ of an alternative algebra $A$ consists
from all quasiregular elements of $A$.

Any simple alternative algebra is Cayley-Dickson algebra $\C(\mu,
\beta, \gamma) = \C$ over their centre $Z(\C)$. As $1 \in Z(\C)$
then $Z(\C) \neq 0$. In such case the centre $Z(\C)$ is a field.
Then $F$ is a field.

Any Cayley-Dickson algebra $\C$ is primitive. The Kleinfeld
radical $\mathcal{K}(A)$ of alternative algebra $A$ is the
intersection of all ideals $I$ of $A$ such that the
quotient-algebra $A/I$ is a primitive algebra. As the algebra
$F(Q)$ is simple  then $(0)$ is the unique maximal ideal of
$F(Q)$. Hence $\mathcal{K}(F(Q)) = (0)$.

The Kleinfeld  radical  $\mathcal{K}(A)$ coincides with Smiley
radical $\mathcal{S}(A)$ in any  alternative algebra $A$,
$\mathcal{K}(A) = \mathcal{S}(A)$. Then from  $\mathcal{K}(F(Q)) =
(0)$ it follows that $\mathcal{S}(F(Q)) = (0)$. But before we
proved that $\mathcal{S}(F(Q)) \neq (0)$. We get a contradiction.
Hence   our supposition that the loop $Q$ is imbedded into a loop
$\mathcal{U}(A)$ for a certain alternative algebra is false. This
completes the proof of Theorem.
\smallskip

\textbf{Remark.} The question on embedding of Moufang loops into
alternative algebras is examined and in Theorem 1 from \cite{San}.
But with regret the statement of these theorem is not correct, is
not in line with the proof. Maybe the translation into English is
not correct. Without getting into details, this theorem is a
consequence of the   following statement that will be published in
the following papers of the author:  the Moufang loops considered
in Theorem of this paper they and only they are not imbedded into
the loop of type $\mathcal{U}(A)$ for a suitable unital
alternative algebra.

\bigskip

\begin{flushleft}
Sandu Nicolae Ion,\\[3mm]

Tiraspol State University of Moldova,\\

Iablochkin str. 5,\\

Kishinev MD-2045, Moldova\\[3mm]

e-mail: {\em sandumn@yahoo.com} \\
\end{flushleft}
\end{document}